 \documentclass[smallabstract,smallcaptions]{dccpaper}
\usepackage{algorithm} 
\usepackage{algpseudocode} 
\usepackage{epsfig}
\usepackage{citesort}
\usepackage{amsmath,amsthm}
\usepackage{amssymb}
\usepackage{color}
\usepackage{url}
\usepackage{soul}
\usepackage{siunitx} 
\newlength{\figurewidth}
\newlength{\smallfigurewidth}
\usepackage{graphicx}
\usepackage{subcaption}
\usepackage{float}
\usepackage{adjustbox}

\usepackage{booktabs} 
\newtheorem{lemma}{Lemma}

\begin{document}

\title
{\large
\textbf{A Modified Levenberg-Marquardt Algorithm For Tensor CP Decomposition in Image Compression }
}

\author{%
Ramin Goudarzi Karim$^{\ast}$, Dipak Dulal$^{\dag}$, and Carmeliza Navasca$^{\ddag}$\\[0.5em]
{\small\begin{minipage}{\linewidth}\begin{center}
\begin{tabular}{ccc}
$^{\ast}$Stillman College & \hspace*{0.5in} & $^{\dag}$ $^{\ddag}$University of Alabama at Birmingham \\
3601 Stillman Blvd && 1720 University Blvd \\
Tuscaloosa, Alabama, 35401, US && Birmingham, Alabama, 35294, US\\
\url{rkarim@stillman.edu} && \url{dpdulal@uab.edu}, \url{cnavasca@uab.edu}
\end{tabular}
\end{center}\end{minipage}}
}

\maketitle
\thispagestyle{empty}

\begin{abstract}
This paper proposed a new variant of the Levenberg-Marquardt algorithm used for Tensor Canonical Polyadic (CP) decomposition with an emphasis on image compression and reconstruction. Tensor computation, especially CP decomposition, holds significant applications in data compression and analysis. In this study, we formulate CP as a nonlinear least squares optimization problem. Then, we present an iterative Levenberg-Marquandt (LM)--based algorithm for computing the CP decomposition. Ultimately, we test the algorithm on various datasets, including randomly generated tensors and RGB images. The proposed method proves to be both efficient and effective, offering a reduced computational burden when compared to the traditional Levenberg-Marquardt technique.
\end{abstract}

\Section{Introduction}

Tensor computation is essential in the domain of computational sciences. Tensors, which can be thought of as higher-order versions of matrices, offer important tools in scientific computing, computer vision, and various engineering disciplines \cite{khoromskij2012tensors,grasedyck2013literature,karim2020accurate}. Furthermore, tensor decompositions constitute essential tools for data science and machine learning applications \cite{kolda2009tensor,sidiropoulos2017tensor}. 

Canonical Polyadic (CP) (a.k.a. CANDECOMP/PARAFAC \cite{hitchcock1927expression}) is a method of tensor factorization. In the simplest terms, tensor factorization involves breaking down a tensor into a product of lower-dimensional tensors. For instance, a third-order tensor can be factorized into a sum of outer products of three-dimensional vectors. The CP decomposition can be useful in many applications, including signal processing \cite{sidiropoulos2017tensor,jiang2022low}, data compression \cite{veganzones2015nonnegative}, and machine learning \cite{zhou2019tensor}, because it allows for a compressed, yet informative, representation of the data. However, finding the best CP decomposition (i.e., minimizing the difference between the original tensor and the decomposed version) can be computationally challenging, particularly for high-dimensional tensors or when the rank is unknown.

The CP decomposition represents a tensor as a summation of component tensors, each with a rank of one. In particular, for a third-order tensor $\mathcal{X} \in \mathbb{R}^{I \times J \times K}$, the CP decomposition allows us to express it as:
\begin{equation} \label{CPproblem}
\mathcal{X} = \sum_{r=1}^{R} a_r \circ b_r \circ c_r
\end{equation}
here, $a_r \in \mathbb{R}^{I}$, $b_r \in \mathbb{R}^{J}$, and $c_r \in \mathbb{R}^{K}$ are vectors, and $\circ$ denotes the outer product \cite{kolda2009tensor}. This decomposition allows us to break down the tensor $\mathcal{X}$ into $R$ simpler rank-one tensors, aiding in the processing of high-dimensional data.
The smallest positive integer $R$ for which the equation (\ref{CPproblem}) holds exactly is referred to as the rank of the tensor \cite{hitchcock1927expression}. Contrary to matrices, the computation of tensor rank poses significant challenges, see \cite{hillar2013most}. Our paper does not tackle these challenges; instead, we focus on computing the CP decomposition of a tensor with a predefined number of rank-one components. The optimization problem related to (\ref{CPproblem}) can be formulated as follows:
\begin{equation}\label{CPopt}
\min_{a_r, b_r, c_r} \quad   \frac{1}{2} \left\Vert \mathcal{X} - \sum_{r=1}^R a_r \circ b_r \circ c_r \right\Vert_F^2 
\end{equation}
where $\Vert \cdot \Vert_F$ represent the Frobenius norm. The optimization problem (\ref{CPopt}) is ill-posed and does not have a unique solution. The ill-posedness of the problem arises from what is referred to as scaling indeterminacy in tensor literature \cite{kolda2009tensor}. Some other authors have employed different formulations equipped with nuclear norms to address this issue. Unlike the formulation in (\ref{CPopt}), the nuclear norm problem leads to a lower semi-continuous function, which possesses a unique solution \cite{friedland2018nuclear,yuan2016tensor}.



The main contributions of this paper are summarized as follows:
\begin{itemize}
    \item We introduce an innovative iteration of the LM algorithm tailored for the CP decomposition. Our approach addresses the computational intensity typically associated with traditional methods by streamlining the process for greater efficiency.
    \item A comprehensive evaluation is conducted using various datasets, ranging from synthetically generated tensors to RGB images.
\end{itemize}
This paper is structured as follows. In Section 1, we review the Levenberg-Marquardt (LM) algorithm and formulate the CP problem as a nonlinear least squares problem. Section 2 describes our proposed algorithm. Section 3 discusses numerical experiments, highlighting the performance and capabilities of our proposed method. Finally, Section 4 concludes the paper, summarizing our key findings and offering insights for potential future work.

\section{Problem Formulation}
The formulation of CP as a nonlinear least squares was first introduced by Paatero \cite{paatero1997weighted}. Subsequent authors have extended this approach to various tensor-related problems \cite{tomasi2005parafac,zhao2023levenberg,chang2023tensor}. In this section, we reformulate the problem (\ref{CPopt}) within the nonlinear least squares framework and derive the essential tools required for the development of our algorithm.\\

\noindent Let $R$  be a fixed positive integer. We introduce the vector-valued function \( F: \mathbb{R}^{R(I+J+K)} \to \mathbb{R}^{IJK} \) characterized by its component functions as:
\begin{equation}\label{defF}
F_{\alpha(i_1,i_2,i_3)}(x) = x_{i_1 i_2 i_3} - \sum_{r=1}^R a_r(i_1) b_r(i_2) c_r(i_3),
\end{equation}
where \( \alpha(i_1, i_2, i_3) = i_1 + (i_2-1)I + (i_3-1)IJ \).
The vector \( x \in \mathbb{R}^{R(I+J+K)} \) is constructed by concatenating the vectors \( a_r \), \( b_r \), and \( c_r \), resulting in:
\begin{equation}
    x = \begin{pmatrix}
        a_1^T &\hdots& a_R^T & b_1^T &\hdots &b_R^T &c_1^T  &\hdots & c_R^T
    \end{pmatrix}^T.
\end{equation}
The factor matrices corresponding to (\ref{CPopt}), expressed in terms of the vectors used to form \( x \), are:
\begin{equation} \label{FacMat}
    A=\begin{pmatrix}
        a_1 & \cdots & a_R
    \end{pmatrix} , \quad B=\begin{pmatrix}
        b_1 & \cdots & b_R
    \end{pmatrix} , \quad C=\begin{pmatrix}
        c_1 & \cdots & c_R
    \end{pmatrix},
\end{equation}
Therefore, the variable \( x \) is defined as the column-wise stacking of the columns of matrices \( A \), \( B \), and \( C \), that is,
\begin{equation}\label{VecX}
x^T = [ \text{vec}(A)^T \text{vec}(B)^T \text{vec}(C)^T] 
\end{equation}
where \( \text{vec}(\cdot) \) denotes the vectorization of a matrix, converting it into a column vector by stacking its columns. The problem (\ref{CPopt}) can now be reformulated as the following nonlinear least squares problem:
\begin{equation} \label{nlsproblem}
    \min _{x \in \mathbb{R}^{R(I+J+K)}}\,  \frac{1}{2} \Vert F(x) \Vert_2^2 .
\end{equation}
The reformulation enables the use of established optimization techniques for nonlinear least squares problems, offering both computational efficiency and insight into the problem's structure. There are numerous techniques available for this specific category of optimization problems \cite{madsen2004methods}.
In this work, we specifically opt for the LM method to address the problem (\ref{nlsproblem}) primarily due to specific characteristics of the Jacobian matrix associated with this problem. The following lemma details the structure of the Jacobian matrix of \( F \) for the problem (\ref{nlsproblem}).
\begin{lemma} \cite{acar2011scalable}\label{lemma:JacobianStructure}
The Jacobian matrix, \( J \), of the function \( F \) in (\ref{nlsproblem}) can be expressed as:
\begin{equation}
J = \begin{pmatrix}
 J_a & J_b & J_c
\end{pmatrix}.
\end{equation}

Each component, \( J_a \), \( J_b \), and \( J_c \), can be further decomposed as:
\begin{align}
J_a &= \begin{pmatrix}
J_a^1 & J_a^2 & \hdots & J_a^R
\end{pmatrix}, 
J_b &= \begin{pmatrix}
J_b^1 & J_b^2 & \hdots & J_b^R
\end{pmatrix}, 
J_c &= \begin{pmatrix}
J_c^1 & J_c^2 & \hdots & J_c^R
\end{pmatrix}.
\end{align}

The individual matrices, \( J_a^r \), \( J_b^r \), and \( J_c^r \), are defined as:
\begin{equation}\label{Jstructure}
J_a^r=-  c_r \otimes b_r \otimes I, \quad
   J_b^r= - c_r\otimes I \otimes a_r, \quad
    J_c^r= -I \otimes  b_r \otimes a_r ,
\end{equation}
\end{lemma}

\begin{figure}[h]
    \centering
    \begin{subfigure}[b]{0.45\textwidth}
        \centering
        \includegraphics[scale=0.45]{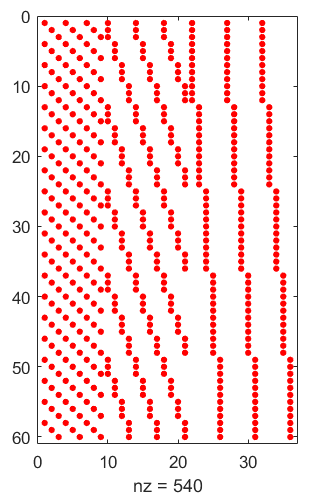}
        \caption{}
        \label{fig:subfig1}
    \end{subfigure}
    \hfill 
    \begin{subfigure}[b]{0.45\textwidth}
        \centering
        \includegraphics[scale=0.45]{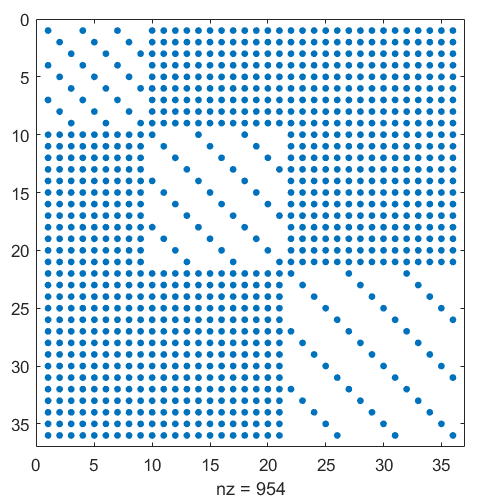}
        \label{fig:subfig2}
        \caption{}
    \end{subfigure}
    \caption{(a) The sparsity of the Jacobian matrix for a \(3 \times 4 \times 5\) tensor with \(R = 2\). The red points indicate the nonzero entries of \(J\). (b) The symmetry and sparsity of the $J^TJ$ for a \(3 \times 4 \times 5\) tensor with \(R = 2\). The blue points indicate the nonzero entries.}
    \label{fig:JJ}
\end{figure}

Here \( \otimes \) represents the Kronecker product of vectors. The Jacobian matrix \( J \) is of size \( Q \times P \), and, due to the appearance of the identity matrix in equation (\ref{Jstructure}), it has a sparse structure. It has $3RQ$ nonzero elements. Figure \ref{fig:JJ} illustrates the Jacobian matrix for a \(3 \times 4 \times 5\) tensor with an estimated rank of \(R=3\). The matrix has dimensions \(60 \times 36\). In total, there are \(3 \times 3 \times 60 = 540\) non-zero elements.


Recall that the LM scheme determines the step size \( h_{LM} \) at each iteration by solving the following linear system of equations:
\begin{equation}\label{LM}
    \left( J^TJ + \mu I \right) h_{LM} = -J^TF
\end{equation}
where \(\mu > 0\) is the \textit{damping parameter}. This parameter plays a crucial role throughout the iterations. It can be interpreted as a regularization parameter. 
The quadratic convergence rate of the LM method, particularly when \(\mu_k = \Vert F_k \Vert^{\delta}\), $\delta \in [1,2]$, has been discussed in \cite{fan2005quadratic,fan2009note}.
Another strategy for determining the damping parameter involves the use of a \textit{gain ratio}. This metric serves as an indicator of the accuracy with which the linear model approximates the objective function. Specifically, if the linear model at the $k$th
  iteration offers a satisfactory approximation, it may be useful to reduce the damping parameter. Conversely, if the approximation is deemed inadequate, an increase in the damping parameter might be warranted \cite{madsen2004methods}.
  In addition to the aforementioned issue, two major computational challenges arise when employing the LM method. Firstly, the dimensionality of the Jacobian matrix can become exceedingly large, particularly for substantial datasets. Specifically, with a size of $IJK \times R(I+J+K)$, the matrix demands significant memory storage at every iteration, posing potential constraints on computational resources. To address this challenge, Tomasi and Bro \cite{tomasi2006comparison}, proposed working directly with \(J^TJ\) and \(J^T F\) rather than evaluating \(J\) at every iteration. By leveraging the symmetric and sparse structure of \(J^TJ\), they optimized computational efforts and significantly reduced memory requirements. Figure \ref{fig:JJ} illustrates the
$J^TJ$ for a $3 \times  4 \times  5$ tensor with an estimated rank of $R = 3$. The matrix
has dimensions $ 36 \times 36$. In total, there are 954 nonzero elements.


The other significant challenge is the high computational requirement of calculating the entries of the Jacobian, as detailed in Lemma \ref{lemma:JacobianStructure}. There are several strategies to mitigate this issue. One approach is to bypass the direct evaluation of the Jacobian matrix and instead approximate it using available quasi-Newton methods. However, it is important to note that adopting this strategy would reduce the convergence rate to superlinear. Another approach, particularly effective when we are sufficiently close to the solution, is to utilize the current Jacobian $J_k$ for the subsequent few iterations, alleviating the need for constant recomputation.

In response to the latter challenge, inspired by the methods in \cite{fan2019adaptive}, we introduce an iterative method utilizing the LM scheme to obtain the CP decomposition of a third-order tensor. This methodology can be naturally extended to tensors of order \( n \). At each iteration \( k \), we use the current Jacobian \( J_k \) to compute both the LM step \( h_k \) and an approximate LM step \( \hat{h}_k \). That is, we solve the following normal equation:
\begin{equation} \label{apph}
 \left( J_k^T J_k + \mu_k I \right) \hat{h}_k = -J_k^T F(y_k), \quad  y_k = x_k + h_k.
\end{equation}


In the subsequent section, we introduce a modified LM algorithm tailored for solving the CP decomposition, which builds upon this approach.

\section{A Modified Levenberg-Marquardt Algorithm for CP Decomposition}

Recall that the Gauss–Newton method \cite{chong2023introduction} is based on the linear approximation of the function \( F \) presented in (\ref{defF}), in the vicinity of \( x_k \). That is,
\begin{equation}\label{LinMod}
    \ell(h) = F(x_k) + J_k h.
\end{equation}
This gives us an approximation of the nonlinear least squares problem (\ref{nlsproblem}) as follows:
\begin{align}
     L(h)&= \frac{1}{2} \ell(h)^T \ell(h) \nonumber\\
     &=\frac{1}{2}F(x_k)^T F(x_k)+h^TJ_k^TF(x_k)+ \frac{1}{2}h^TJ_k^TJ_kh
\end{align}
If \( J_k \) is full rank, then the unique minimizer of \( L(h) \) can be found by solving the following normal equation:
\begin{equation}
    \left( J_k^T J_k \right) h = - J_k^T F(x_k)
\end{equation}
However, the Jacobian matrix in Lemma \ref{lemma:JacobianStructure} does not possess full rank. Specifically, the rank of this matrix is always bounded above by
$P-2R$. This issue can be addressed by introducing a damping parameter $\mu_k $ in each iteration. The damped Gauss-Newton (dGN) algorithm addresses this by solving the following normal equation in each iteration:
\begin{equation} \label{LMStep}
    \left( J_k^T J_k +\mu_k I\right) h = - J_k^T F(x_k).
\end{equation}
The inclusion of the damping parameter \( \mu_k \) in equation (\ref{LMStep}) ensures that the matrix \( J_k^T J_k + \mu_k I \) is positive definite, which consequently guarantees that the step size \( h \) is uniquely determined for each iteration. This approach can also be interpreted within the framework of trust region methods. Specifically, we are effectively solving the following constrained optimization problem at each iteration:

\begin{equation}
    \min_h \, L(h) \quad \text{s.t.} \quad \Vert h \Vert_2 \leq \Delta_k,
\end{equation}
where \( \Delta_k \) is the trust region radius, analogous to \( \mu_k \), that determines the size of the region within which the model is trusted to be an accurate representation of the objective function.
In order to make use of the current Jacobian $J_k$ we also predict the next step $\hat{h}_k$
by solving the normal equation (\ref{apph}) and set the new step size $s_k=h_k+\hat{h_k}$. The trail step $y_{k+1}=x_k+s_k$ is considered. 
In order to make use of the current Jacobian \( J_k \), we also predict the next step \( \hat{h}_k \)
by solving the normal equation (\ref{apph}) and set the new step size \( s_k = h_k + \hat{h}_k \). The trial step \( y_{k+1} = x_k + s_k \) is then considered.
Note that in equation \eqref{apph}, there is no need to compute the Jacobian matrix at \( y_k \). Only the function value \( F(y_k) \) is required to solve \eqref{apph}. As a result, the cost of obtaining \( \hat{h}_k \) will be relatively inexpensive. This can be especially advantageous when the current iteration \( k \) is close enough to the solution.
 To test whether or not to accept the trial step, we introduce the gain ratio at every iteration:
\begin{equation} \label{GainRat}
    \rho_k = \frac{\text{actual reduction}}{\text{predicted reduction}} =\frac{\Vert F(x_k)\Vert - \Vert F(x_k+s_k) \Vert}{\Vert F(x_k) \Vert - \Vert \ell (h_k) \Vert + \Vert F(y_k) \Vert - \Vert \ell(\hat{h}_k) \Vert}
\end{equation}
If \(\rho_k\) is greater than a predefined threshold or is at least positive, we accept the trial step as the linear model (\ref{LinMod}) is working satisfactorily. Conversely, if \(\rho_k\) is less than the gain threshold, the linear model is not a good approximation of \( F \), so the trial step is rejected. Consequently, we increase the damping parameter \(\mu\) by multiplying it with a constant factor.
Algorithm \ref{alg:MLM} summarizes the method described in this section.

\begin{table}[htb]
  \centering
  
  \label{tab:comparison}
  \begin{tabular}{
    |l|
    c|
    S[table-format=2.2]|
    S[table-format=1.2e-1]|
    S[table-format=2.2]|
    S[table-format=1.2e-1]|
    S[table-format=1.2]|
  }
    \hline
    \textbf{tensor size} & \textbf{rank} & \multicolumn{2}{c|}{\textbf{LM}} & \multicolumn{2}{c|}{\textbf{Modified LM}} & \textbf{comp} \\
    \cline{3-6}
    & & \textbf{time } & \textbf{residual error} & \textbf{time} & \textbf{residual error} &\\
    \hline
   $35\times 25\times 15$  & 40  & 65.8 &306.4973  & 49.3 & 306.497 & 77 \\
    \hline
    $20 \times 20 \times 12$ & 30 & 11.4 & 84.706& 7.4 &  84.710& 68 \\
    \hline
   $28 \times 18 \times 16$ & 35  & 27 & 163.83 & 19.1 & 163.94 & 73 \\
    \hline
  \end{tabular}
  \caption{LM vs Modified LM }
  \label{tab:3}
\end{table}

\begin{table}[htb]
    \centering

    \begin{tabular}{|c|c|c|c|c|c|}
        \toprule
        \textbf{image}  & \textbf{rank} & \textbf{time(m)} & \textbf{residual error} & \textbf{comp }\\
        \midrule
        Rose & 20\textbackslash50\textbackslash75 & 3.4\textbackslash14.5\textbackslash25.9 & 157.4\textbackslash13.1\textbackslash2.07 & 87\textbackslash67 \textbackslash50 \\
        \midrule
        Leopard & 20\textbackslash50\textbackslash80 & 19.5\textbackslash78\textbackslash178 & 165\textbackslash83.7\textbackslash53 & 92\textbackslash79\textbackslash67 \\
        \midrule
        Pepper & 25\textbackslash50\textbackslash80 & 30.1\textbackslash95.2\textbackslash209 & 302.5\textbackslash302.5\textbackslash95.2\textbackslash29.3 & 90\textbackslash79\textbackslash68 \\
        \bottomrule
    \end{tabular}
    \caption{\textbf{ Modified LM:} time vs. error vs. compression on RGB images}
    \label{tab:1}
\end{table}


\begin{footnotesize}
\begin{algorithm}[H]
	\caption{Modified LM for CP decomposition} \label{alg:MLM}
	\begin{algorithmic}[1]\scriptsize
    \State \textbf{input:} Initial factor matrices $A$, $B$ and $C$, third order tensor $\mathcal{X}$
    \State \textbf{stopping criteria:} max iteration $N$, error tolerance $\epsilon$
    \State \textbf{parameters:} damping parameter $\mu>0$, gain ratio constant $\gamma$, damping parameter multiple $\nu>1$
    \State \textbf{initialize:} evaluate the initial $x$ and $J$ by equations (\ref{VecX}) and (\ref{lemma:JacobianStructure})
		\While {stopping criteria not met}
		\State \textbf{solve:} \[\left( J_k^T J_k +\mu I\right) h_k = - J_k^T F(x_k)\]
  \State \textbf{set:}
  \[ y_k = x_k + h_k\]
        \State \textbf{solve:}
        \[\left( J_k^T J_k + \mu I \right) \hat{h}_k = -J_k^T F(y_k), \quad \]

        \State \textbf{set the trial step:}
        \[s_k=h_k+\hat{h}_k, \quad y_{k+1} =x_k+s_k\]
        \State \textbf{compute:} $\rho_k$ as in (\ref{GainRat})
        \If {$\rho_k > \gamma$}
        \State accept trial step and decrease damping parameter:
        \[x_{k+1}=y_{k+1},\quad  \text{reset} \, \nu, \quad \mu=(1/2)\mu \]
        \Else
        \State reject trial step and increase damping parameter:
        \[\mu= \nu \mu, \quad \nu=2\nu\]
        \EndIf
  \EndWhile
  \For {$r=1:R$}
          \State  \[A_{est}(:,r)=x((r-1)I+1:rI)\]
            \[B_{est}(:,r)=x(RI+(r-1)J+1:RI(r)J)\]
            \[C_{est}(:,r)=x(R(I+J)+(r-1)K+1:R(I+J)+(r)K)\]    
    
  \EndFor
  \State \textbf{output:} three factor matrices $A_{est}, B_{est}$ and $C_{est}$ of the CP decomposition.
	\end{algorithmic} 
\end{algorithm} 
\end{footnotesize}

\section{Numerical Experiments}

In this section, we evaluate and compare the performance of the proposed algorithm through a set of numerical experiments from the standard tensor-based Levenberg-Marquandt method. The computations were carried out on UAB's high-performance computing system (Cheaha), running MATLAB 2023a with a Pascalnodes node partition and 18 GB of memory per CPU with 2 GPUs and 12 CPUs per node. In our numerical experiments, we tested two types of data: RGB images and randomly generated tensors with varying rank, mode dimensions, and other features. 

\begin{figure}[H]
    \centering
    \begin{subfigure}[b]{0.45\textwidth}
        \centering
        \includegraphics[scale=0.45]{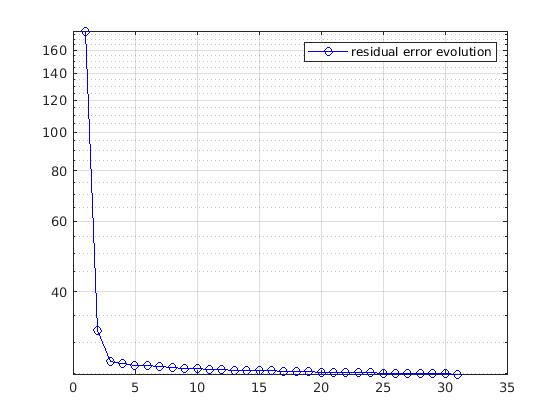}
        \caption{}
        \label{fig:subres1}
    \end{subfigure}
    \hfill 
    \begin{subfigure}[b]{0.45\textwidth}
        \centering
        \includegraphics[scale=0.45]{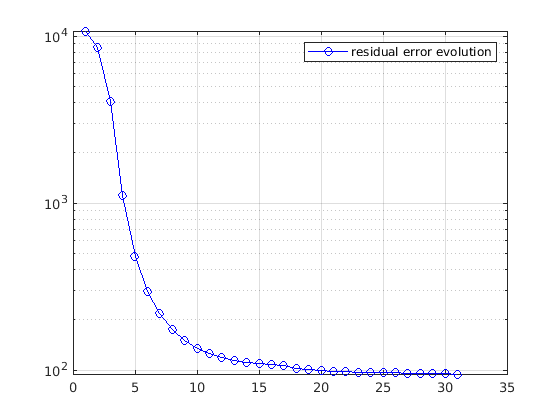}
        \label{fig:subres2}
        \caption{}
    \end{subfigure}
    \caption{\textbf{Modified LM}: The performance of the proposed algorithm in terms of accuracy of the residual error for (a) a randomly generated tensor ($45\times35\times 25),$ R=40, 89\% compression and (b) Leopard image R = 50, 79\%  compression.}
    \label{residual_er}
\end{figure}

\begin{figure}[H]
  \begin{subfigure}{0.3\textwidth}
    \includegraphics[width = \linewidth]{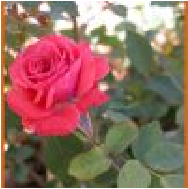}
    \caption{original}
    \label{subfig1}
  \end{subfigure}
  \hfill
  \begin{subfigure}{0.3\textwidth}
    \includegraphics[width = \linewidth]{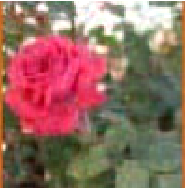}
    \caption{$87\%$ ,R = 20}
    \label{subfig}
  \end{subfigure}
  \hfill
  \begin{subfigure}{0.3\textwidth}
    \includegraphics[width = \linewidth]{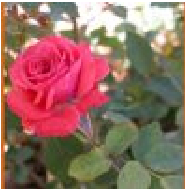}
    \caption{$50\%$, R = 75}
  \end{subfigure}
\hfill
  \begin{subfigure}{0.3\textwidth}
    \includegraphics[width=\linewidth]{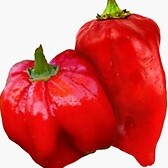}
    \caption{original}
  \end{subfigure}
  \hfill
  \begin{subfigure}{0.3\textwidth}
    \includegraphics[width=\linewidth]{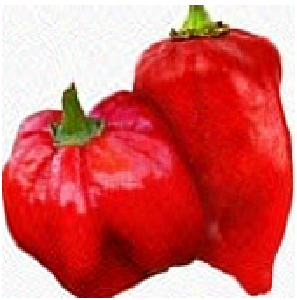}
    \caption{$80\%$,R = 50}
  \end{subfigure}
  \hfill
  \begin{subfigure}{0.3\textwidth}
    \includegraphics[width=\linewidth]{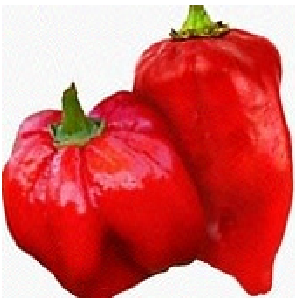}
    \caption{$68\%$,R = 80}
  \end{subfigure}
  \hfill
   \begin{subfigure}{0.3\textwidth}
    \includegraphics[width=\linewidth]{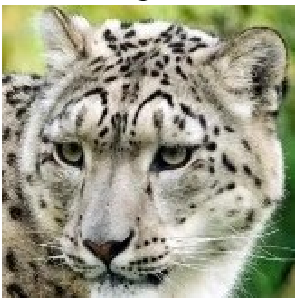}
    \caption{original}
  \end{subfigure}
  \hfill
  \begin{subfigure}{0.3\textwidth}
    \includegraphics[width=\linewidth]{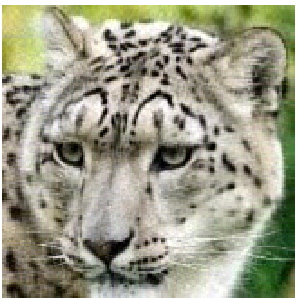}
    \caption{$79\%$ R= 50}
  \end{subfigure}
  \hfill
  \begin{subfigure}{0.3\textwidth}
    \includegraphics[width=\linewidth]{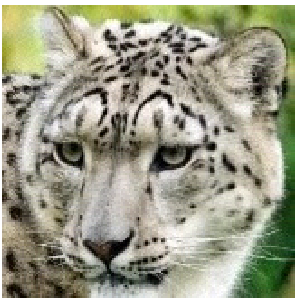}
    \caption{$67\%$, R = 80}
  \end{subfigure}
  \caption{\textbf{Modified LM}: Rose($100\times 100\times 3$), Pepper ($168\times 168\times 3$) and Leopard($162\times 162\times 3$).}
  \label{rgb}
\end{figure}

To demonstrate the effectiveness of the Modified LM method, various third-order tensors of RGB images were reconstructed. 
In Figure \ref{rgb}, the Modified LM reconstructed three RGB images: rose, pepper, and leopard. The compression percentage is calculated from the ratio: $\frac{R(I + J + K)}{I \times J \times K}$ for a third order tensor of size $I \times J \times K$ of rank $R$. The ratio measures how much an image has been reduced in size. Clearly, the higher the compression percentage, the greater the reduction in the size of the image. In Figure \ref{rgb}, the images in the middle and last columns are reconstructed tensor images with varying ranks and, hence, compression percentages. The images look indistinguishable from the original images. Table \ref{tab:1} shows that a higher compression requires less time but with low accuracy. Similarly, Table \ref{tab:2} follows the same trend in Table \ref{tab:1} with randomly generated tensors for some specific sizes. Lastly, Table \ref{tab:3} compares the two algorithms, LM and modified LM, for random tensors of sizes: $35 \times 25 \times 15$ with $R=40$, $20 \times 20 \times 12$ with $R=30$, and $28 \times 18 \times 16$ with $R=35$. In all three cases, CPU times for modified LM is much less than for LM. The smooth evolution of residual errors during the model execution is shown in Figure \ref{residual_er}.
\vspace{1cm}

\begin{table}[H]
\centering
\begin{tabular}{@{}|c|c|c|c|c|c|@{}}
\toprule
\textbf{tensor size} & \textbf{rank} & \textbf{time cost(m)} & \textbf{residual error} & \textbf{comp}\\
\midrule
45 $\times$ 35 $\times$ 20 & 25\textbackslash40\textbackslash60 & 1.9\textbackslash3.5\textbackslash5.9 & 21.7\textbackslash19.5\textbackslash16.7 & 91\textbackslash87\textbackslash81 \\
\midrule
45 $\times$ 35 $\times$ 25 & 25\textbackslash40\textbackslash60  & 2.6\textbackslash4.9\textbackslash7.5  & 24.8\textbackslash22.8\textbackslash20.08 & 93\textbackslash89\textbackslash84  \\
\midrule
45 $\times$ 35 $\times$ 30 & 25\textbackslash40\textbackslash60& 3\textbackslash5.4\textbackslash9.1 & 27.7\textbackslash25.6\textbackslash22.0 & 94\textbackslash91\textbackslash86 \\
\bottomrule
\end{tabular}
\caption{\textbf{Modified LM}: time vs. error vs. compression on randomly generated third oder tensors}
\label{tab:2}
\end{table}

\section{Conclusion}
In this work, we proposed a new method for reconstructing higher-order low rank tensors in the canonical polyadic decomposition format. We build on the success of the Levenberg-Marquardt approach on approximating a well-scaled solution of a nonlinear least-squares problem. However, there are drawbacks and challenges in the LM method. Our approach is a modified LM which decreases the computational demand on the calculation of the Jacobian and objective function evaluation. From numerical experiments, the modified LM performed well in compressing RGB images and randomly generated tensors while keeping the tensor rank low.  The modified LM has significantly lower CPU times in all experimental cases compared to the LM. In conclusion, a modified LM provides sensible compression for low-rank CP reconstruction of tensors with the same order of accuracy as the LM method. This is achieved while leveraging a noticeable reduction in CPU time cost. In our future outlook, we will add sampling methods to tackle massively high mode dimensions in high-order tensors and will provide theoretical convergence results.

\Section{References}
\bibliographystyle{IEEEbib}
\bibliography{ref}

\begin{thebibliography}{10}

\bibitem{khoromskij2012tensors}
Boris~N Khoromskij,
\newblock ``Tensors-structured numerical methods in scientific computing: Survey on recent advances,''
\newblock {\em Chemometrics and Intelligent Laboratory Systems}, vol. 110, no. 1, pp. 1--19, 2012.

\bibitem{grasedyck2013literature}
Lars Grasedyck, Daniel Kressner, and Christine Tobler,
\newblock ``A literature survey of low-rank tensor approximation techniques,''
\newblock {\em GAMM-Mitteilungen}, vol. 36, no. 1, pp. 53--78, 2013.

\bibitem{karim2020accurate}
Ramin~Goudarzi Karim, Guimu Guo, Da~Yan, and Carmeliza Navasca,
\newblock ``Accurate tensor decomposition with simultaneous rank approximation for surveillance videos,''
\newblock in {\em 2020 54th Asilomar Conference on Signals, Systems, and Computers}. IEEE, 2020, pp. 842--846.

\bibitem{kolda2009tensor}
Tamara~G Kolda and Brett~W Bader,
\newblock ``Tensor decompositions and applications,''
\newblock {\em SIAM review}, vol. 51, no. 3, pp. 455--500, 2009.

\bibitem{sidiropoulos2017tensor}
Nicholas~D Sidiropoulos, Lieven De~Lathauwer, Xiao Fu, Kejun Huang, Evangelos~E Papalexakis, and Christos Faloutsos,
\newblock ``Tensor decomposition for signal processing and machine learning,''
\newblock {\em IEEE Transactions on signal processing}, vol. 65, no. 13, pp. 3551--3582, 2017.

\bibitem{hitchcock1927expression}
Frank~L Hitchcock,
\newblock ``The expression of a tensor or a polyadic as a sum of products,''
\newblock {\em Journal of Mathematics and Physics}, vol. 6, no. 1-4, pp. 164--189, 1927.

\bibitem{jiang2022low}
Jiahua Jiang, Fatoumata Sanogo, and Carmeliza Navasca,
\newblock ``Low-cp-rank tensor completion via practical regularization,''
\newblock {\em Journal of Scientific Computing}, vol. 91, no. 1, pp. 18, 2022.

\bibitem{veganzones2015nonnegative}
Miguel~A Veganzones, Jeremy~E Cohen, Rodrigo~Cabral Farias, Jocelyn Chanussot, and Pierre Comon,
\newblock ``Nonnegative tensor cp decomposition of hyperspectral data,''
\newblock {\em IEEE Transactions on Geoscience and Remote Sensing}, vol. 54, no. 5, pp. 2577--2588, 2015.

\bibitem{zhou2019tensor}
Mingyi Zhou, Yipeng Liu, Zhen Long, Longxi Chen, and Ce~Zhu,
\newblock ``Tensor rank learning in cp decomposition via convolutional neural network,''
\newblock {\em Signal Processing: Image Communication}, vol. 73, pp. 12--21, 2019.

\bibitem{hillar2013most}
Christopher~J Hillar and Lek-Heng Lim,
\newblock ``Most tensor problems are np-hard,''
\newblock {\em Journal of the ACM (JACM)}, vol. 60, no. 6, pp. 1--39, 2013.

\bibitem{friedland2018nuclear}
Shmuel Friedland and Lek-Heng Lim,
\newblock ``Nuclear norm of higher-order tensors,''
\newblock {\em Mathematics of Computation}, vol. 87, no. 311, pp. 1255--1281, 2018.

\bibitem{yuan2016tensor}
Ming Yuan and Cun-Hui Zhang,
\newblock ``On tensor completion via nuclear norm minimization,''
\newblock {\em Foundations of Computational Mathematics}, vol. 16, no. 4, pp. 1031--1068, 2016.

\bibitem{paatero1997weighted}
Pentti Paatero,
\newblock ``A weighted non-negative least squares algorithm for three-way ‘parafac’factor analysis,''
\newblock {\em Chemometrics and Intelligent Laboratory Systems}, vol. 38, no. 2, pp. 223--242, 1997.

\bibitem{tomasi2005parafac}
Giorgio Tomasi and Rasmus Bro,
\newblock ``Parafac and missing values,''
\newblock {\em Chemometrics and Intelligent Laboratory Systems}, vol. 75, no. 2, pp. 163--180, 2005.

\bibitem{zhao2023levenberg}
Jinyao Zhao, Xuejuan Zhang, and Jinling Zhao,
\newblock ``A levenberg-marquardt method for tensor approximation,''
\newblock {\em Symmetry}, vol. 15, no. 3, pp. 694, 2023.

\bibitem{chang2023tensor}
Shih~Yu Chang, Hsiao-Chun Wu, Yen-Cheng Kuan, and Yiyan Wu,
\newblock ``Tensor levenberg-marquardt algorithm for multi-relational traffic prediction,''
\newblock {\em IEEE Transactions on Vehicular Technology}, 2023.

\bibitem{madsen2004methods}
Kaj Madsen, Hans~Bruun Nielsen, and Ole Tingleff,
\newblock ``Methods for non-linear least squares problems,''
\newblock {\em lecture notes}, 2004.

\bibitem{acar2011scalable}
Evrim Acar, Daniel~M Dunlavy, and Tamara~G Kolda,
\newblock ``A scalable optimization approach for fitting canonical tensor decompositions,''
\newblock {\em Journal of chemometrics}, vol. 25, no. 2, pp. 67--86, 2011.

\bibitem{fan2005quadratic}
Jin-yan Fan and Ya-xiang Yuan,
\newblock ``On the quadratic convergence of the levenberg-marquardt method without nonsingularity assumption,''
\newblock {\em Computing}, vol. 74, pp. 23--39, 2005.

\bibitem{fan2009note}
Jinyan Fan and Jianyu Pan,
\newblock ``A note on the levenberg--marquardt parameter,''
\newblock {\em Applied Mathematics and Computation}, vol. 207, no. 2, pp. 351--359, 2009.

\bibitem{tomasi2006comparison}
Giorgio Tomasi and Rasmus Bro,
\newblock ``A comparison of algorithms for fitting the parafac model,''
\newblock {\em Computational Statistics \& Data Analysis}, vol. 50, no. 7, pp. 1700--1734, 2006.

\bibitem{fan2019adaptive}
Jinyan Fan, Jianchao Huang, and Jianyu Pan,
\newblock ``An adaptive multi-step levenberg--marquardt method,''
\newblock {\em Journal of Scientific Computing}, vol. 78, pp. 531--548, 2019.

\bibitem{chong2023introduction}
Edwin~KP Chong, Wu-Sheng Lu, and Stanislaw~H {\.Z}ak,
\newblock {\em An introduction to optimization},
\newblock John Wiley \& Sons, 2023.

\end{thebibliography}
\end{document}